\documentstyle[12pt,amscd]{amsart}
\newtheorem{Theorem}{Theorem}[section]
\newtheorem{Lemma}[Theorem]{Lemma}
\newtheorem{Corollary}[Theorem]{Corollary}
\newtheorem{Proposition}[Theorem]{Proposition}

\def\depth{\operatorname{depth}}
\def\reg{\operatorname{reg}}

\def\Ext{\operatorname{Ext}}

\def\geom{\operatorname{g-reg}}
\def\hdeg{\operatorname{hdeg}}

\def\To {\longrightarrow}
\def\sk{\smallskip\par}
\def\mm{{\frak m}}
\def\qq{{\frak q}}

\def\PP{{\Bbb P}}
\def\F{{\cal F}}
\def\M{{\cal M}}
\begin{document}

\title{Castelnuovo-Mumford regularity\\  and extended degree}
\author{Maria Evelina Rossi, Ng\^o Vi\^et Trung and Giuseppe Valla}
\address{Dipartimento di Matematica,  Universit\`a di Genova, Via
Dodecaneso 35, 16132 Genova, Italy}
\email{rossim@@dima.unige.it}
\address{Institute of Mathematics,  Box 631, B\`o H\^o, 10000 Hanoi, Vietnam}
\email{nvtrung@@thevinh.ncst.ac.vn}
\address{Dipartimento di Matematica,  Universit\`a di Genova, Via
Dodecaneso 35, 16132 Genova, Italy}
\email{valla@@dima.unige.it}
\thanks{The first and third authors are partially supported by MPI of
Italy. The second author is partially supported
by the National Basic Research Program of Vietnam}
\begin{abstract} The main result of this paper shows that the
Castelnuovo-Mumford regularity of the tangent cone of a local ring
$A$ is effectively bounded by the dimension and any extended
degree of $A$. From this it follows that there are only a finite
number of Hilbert-Samuel functions of local rings with given
dimension and extended degree.  \end{abstract}
\maketitle

\section*{Introduction} \sk

The Castelnuovo-Mumford regularity is a kind of universal bound for
important invariants of graded algebras such as the maximum degree of
the syzygies  and the maximum non-vanishing degree of the local
cohomology modules (see Section 1).
It has been used as a measure for the complexity of computational
problems in algebraic geometry and commutative algebra (see e.g.
[EG], [BM], [V2]).
One has often tried to find upper bounds for the
Castelnuovo-Mumford regularity in terms of simpler invariants.
The simplest invariants which reflect the complexity of a graded
algebra are the
dimension and the multiplicity. However, the Castelnuovo-Mumford
regularity can not be bounded in
terms of the multiplicity and the dimension. \sk

Extended degree was recently introduced by  Vasconcelos et al [DGV], [V1],
[V2]  in order to capture the size of a module along with some of the
complexity of its structure. It is a numerical function on the
category of finitely generated
modules over local or graded rings which generalizes the usual notion
of multiplicity and degree (see  Section 2).
These invariants  tend to have a homological character but are still
amenable to explicit computation by computer algebra systems.
It has been shown by Doering, Gunston and Vasconcelos [DGV] that the
Castelnuovo-Mumford regularity of a graded algebra is less than any
extended degree. \sk

In this paper we show that the Castelnuovo-Mumford regularity of the
tangent cone of a local  ring is bounded above by an exponential
function of the dimension and any extended degree. This
bound is neither  a consequence nor a simple generalization of the
result of [DGV] since there is no relationship between the extended
degrees
of a local ring and those of its tangent cone. \sk

As an application we give upper bounds for the coefficients of the
Hilbert-Samuel function in terms of any extended degree.
It follows that there are only a finite number of
Hilbert-Samuel functions for local rings with given dimension and
extended degree. This application covers three recent results on the
finiteness of Hilbert functions. The first result is due to Srinivas
and Trivedi  [ST2]  who showed that there exist only a finite number
of Hilbert functions  of
Cohen-Macaulay local rings with given dimension and multiplicity (see
also [KL] [ST1], [Tri1]).
This result was extended by Trivedi [Tri2] to the class of
generalized Cohen-Macaulay rings by involving the lengths of local
cohomology modules. The last result is due to Rossi, Valla and
Vasconcelos [RVV] who
showed that there exist only a finite number of Hilbert functions of
graded algebras with given dimension and extended degree.  \sk

Let $(A,\mm)$ be a local ring with infinite residue field and $d = \dim A$.
We will denote by $\reg(G)$ the Castelnuovo-Mumford regularity of the
associated graded ring
$G = \oplus_{n \ge 0}\mm^n/\mm^{n+1}$ and we will write  the
Hilbert-Samuel polynomial of $A$ in the form
$$P_A(X):=\sum_{i=0}^d(-1)^ie_i(A){X+d-i \choose d-i}.$$
By definition $P_A(n)= \ell(A/\mm^{n+1}) $ for $n \gg 0.$
Moreover, we set $e(A) = e_0(A)$ and $I(A) = D(A)-e(A)$, where $D(A)$
is a given extended degree of $A$. We always have $I(A) \ge 0$ with
$I(A) = 0$ if and only if $A$ is a Cohen-Macaulay ring.   \sk

Our main results  can be formulated as follows.\medskip

\noindent{\bf Theorem \ref{regularity}.}
Let $D(A)$ be an arbitrary
extended degree of $A$. Then\par
{\rm (i) } $\reg(G) \le e(A) + I(A)-1$ if $d = 1$,\par
{\rm (ii)} $\reg(G) \le
   e(A)^{(d-1)!-1}[e(A)^2+e(A)I(A)+2I(A)-e(A)]^{(d-1)!}- I(A)$ if $d \ge 2$.
   \medskip

\noindent{\bf Theorem \ref{coeff}.}  Let $D(A)$ be an arbitrary
extended degree of $A$. Then \par
{\rm (i)} $|e_1(A)| \le  \displaystyle \frac{e(A)[e(A)-1]}{2}+I(A)$,\par
{\rm (ii)} $|e_i(A)| \le e(A)^{i!-i}[e(A)^2+e(A)I(A)+2I(A)]^{i!}-1$ if $i
\ge 2$.
   \medskip

Though these bounds are far from being optimal, they provide the
means to make {\it a priori} estimates in local algebra, following
the philosophy of [DGV] and [V2].\sk

The starting point of our approach is the observation that the
regularity of the tangent cone can be estimated by means of the
geometric regularity (Proposition \ref{depthzero} and Propositon
\ref{Hoa}).
Due to a result of Mumford [M],  the problem of bounding the
geometric regularity can be reduced to  the problem of bounding the
Hilbert polynomial (Theorem \ref{Mumford}). This idea has been used
in [ST2] and [Tri2]. But unlike [ST2] and [Tri2] which refer to deep
results from algebraic geometry such as Grothendieck's formal
function theorem, we use standard  algebraic methods to solve this
problem.  The key point is a uniform bound for  the  Hilbert-Samuel
function in terms of any extended degree (Theorem \ref{uniform}).  By
induction, this bound allows us to estimate the regularity of the
tangent cone. The coefficients of the Hilbert-Samuel polynomial can
be bounded then by a device due to Vasconcelos. \sk

We would like to mention that Trivedi [Tri1], [Tri2]  has dealt with
the Hilbert function of Cohen-Macaulay and generalized Cohen-Macaulay
modules with respect to $\mm$-primary ideals. Our approach shows that
the results of this paper can be also extended to this general
situation.  But it is not our intention to go so far. \sk

The paper is organized as follows. In Section 1 we prepare some facts
and results on the regularity of graded algebras. In Section 2
we estimate the geometric regularity in terms of any extended degree.
In Section 3 and Section 4 we prove the bound for the regularity of
the tangent cone and the bounds for the coefficients of the
Hilbert-Samuel polynomial, respectively.

\section{Regularity of graded algebras}

Throughout this section let $S = k[x_1,\ldots,x_r]$ be a polynomial
ring over a field $k$.
Let $M$ be a finitely generated graded $S$-module. Let
$$0 \To  F_s \To  \cdots \To  F_1 \To  F_0 \To  M \To 0$$
be a minimal graded free resolution of $M$.
Write  $b_i$ for the maximum of the  degrees of the generators of $F_i$.
Following [E, Section 20.5] we say that $M$ is {\it m-regular} for
some integer $m$ if $b_j-j\le
m$ for all $j$.
The {\it Castelnuovo-Mumford regularity} $\reg(M)$ of $M$ is defined
to be the least integer $m$ for which $M$ is $m$-regular, that is,
$$\reg(M) = \max\{b_i-i|\ i = 0,\ldots,s\}.$$

It is well known that $M$ is $m$-regular if and only if
$\Ext_S^i(M,S)_n=0$ for all $i$ and all $n\le -m-i-1$ (see [EG] and
[E, Proposition 20.16]).
This result is hard to apply because in principle infinitely many
conditions must be checked. However, in some cases, it suffices to
check just a few. \sk

We say that $M$ is {\it weakly m-regular} if $\Ext_S^i(M,S)_{-m-i-1}=0$ for
all $i$ [E, Section 20.5]. Concerning
this weaker notion of regularity, we have the following result of Mumford. \sk

\begin{Theorem} \label{weakly} {\rm (see [E, Theorem 20.17])} If $M$
is weakly $m$-regular and $L$ is the
maximal submodule of $M$ having finite length, then $M/L$ is $m$-regular.
\end{Theorem}

If $M$ has positive depth, then $L=0$ so that $m$-regularity and weak
$m$-regularity coincide. It is less known that this
is also the case for homogeneous quotient rings of $S$ (not
necessarily of positive depth). \sk

\begin{Corollary} \label{basic1} Let $R=S/I$, where $I$ is  an
homogeneous ideal. If $R$ is weakly $m$-regular, then $R$ is $m$-regular.
\end{Corollary}

\begin{pf}  By Theorem \ref{weakly},  $R/L$ is $m$-regular,  where
$L$ is the largest ideal of $R$ of finite length. Hence $m \ge 0$
by the definition of regularity. From the short exact sequence
$$0\To  I \To  S \To  R \To  0$$
we get
\begin{eqnarray*} & \Ext_S^0(I,S)_{-m-1} \simeq
\Ext_S^{1}(R,S)_{-m-1} = 0, \\ & \Ext_S^i(I,S)_{-m-i-1} \simeq
\Ext_S^{i+1}(R,S)_{-m-i-1} = 0,\ i \ge 1, \end{eqnarray*}
so that $I$ is weakly $(m+1)$-regular. Since $I$ has positive depth,
$I$ is $(m+1)$-regular by Theorem \ref{weakly}. Looking at the
minimal graded free resolution of $I$, we can conclude that $R$ is
$m$-regular.
\end{pf}

Let $R$ be a standard graded algebra over a field $k$. Let $M$ now be
a finitely generated graded $R$-module.
For any integer $i$ we denote by $H_{R_+}^i(M)$ the $i$-th local
cohomology module of $M$, where $R_+$ is the maximal graded ideal of $R$.
\par

If $R = S/I$, where $I$ is a homogeneous ideal, then $H_{R_+}^i(M) =
H_{S_+}^i(M)$. By local duality (see [E,
A4.2]) we have
$$H_{S_+}^i(M)_m \cong \Ext_S^{r-i}(M,S)_{-m-r}$$ for all $i$ and $m$.
Thus, $M$ is $m$-regular if and only if
$H_{R_+}^i(M)_n =0$ for all $i$ and $n \ge m-i+1$, and
$M$ is weakly $m$-regular if and only if $H_{R_+}^i(M)_{m-i+1}=0$ for all
$i$. In particular,
$\reg(M)$ is the least integer $m$ for which
$H_{R_+}^i(M)_n = 0$ for all $i$ and $n \ge m-i+1$. Hence the
Castelnuovo-Mumford regularity can be defined for any
finite $R$-module regardless of its presentation. \sk

We can use the Castelnuovo-Mumford regularity to control the behavior
of the Hilbert function. Let $h_M(n)$ the
Hilbert function of $M$ and
$p_M(n) $ the Hilbert polynomial of $M$ which is by definition the
unique polynomial $\in {\bf Q}[X] $ for which
$p_M(n)=h_M(n) $ for all $n \gg 0.$ Their difference is determined by
the following classical formula of Serre (see e.g. [BH,
Theorem 4.4.3]):
$$h_M(n) - p_M(n) = \sum_{i\ge 0} (-1)^i \dim_kH_{R_+}^i(M)_n.$$
Hence we immediately obtain the following consequence.

\begin{Lemma} \label{Serre} $h_R(n) = p_R(n)$ for $n> \reg(R)$.
\end{Lemma}

The regularity in algebraic geometry is defined a bit differently.
Let $\F$ be a coherent sheaf on $\PP^r$. Then $\F$ is called
$m$-regular if
$H^i(\PP^r,\F(m-i)) = 0$ for all $i > 0$, where $H^i(\PP^r,\F)$ denotes
the $i$th sheaf cohomology of $\F$ [Mu, p. 99, Definition]. This
regularity is related to the
Castelnuovo-Mumford regularity by the Serre-Grothendieck
correspondence which says that if $\F$ is the sheaf
associated to the $R$-module $M$, then  $H^i(\PP^r,\F(n)) \cong
H_{R_+}^{i+1}(M)_n$ for all $n$ and $i
\ge 1$. \sk

Due to a result of Mumford in [Mu], the regularity of an ideal sheaf
can be estimated in terms of that of a
generic hyperplane section by means of the Hilbert polynomial. We
present here an algebraic version of this
result, for which we shall need the following notations.\sk

\noindent{\bf Definition.} We say that $M$ is {\it geometrically
$m$-regular} if $H_{R_+}^{i}(M)_n=0$ for all $i>0$
and $n \ge m-i+1$, and we define the {\it geometric regularity} $\geom(M)$
of $M$ to be the least integer $m$ for
which $M$ is geometrically $m$-regular.\sk

We always have $\geom(M) \le \reg(M)$ and $\geom(M) = \reg(M)$ if
$\depth M > 0$.

Following [Tru2] we call a homogeneous element $z$ of $R$ {\it
filter-regular} if $z \not\in Q$ for any relevant associated prime
ideal $Q$ of $R$. Filter-regular linear forms always exist if $k$ is an
infinite field.  \sk

The following theorem follows  from the proof of [Mu, pp. 101,
Theorem]. We insert here a proof for
completeness. Note that $R = S/I$ is geometrically $n$-regular if and only
if the sheafification of $I$ is $(n+1)$-regular.  \sk

\begin{Theorem} \label{Mumford} Let $R$ be a standard graded algebra
with $\dim R \ge 1$. Let $z$ be a
filter-regular linear form of $R$. If $R/zR$ is geometrically $m$-regular, then $R$ is geometrically $(m+p_R(m)-h_{R/L}(m))$-regular,
where $L$ denotes the largest ideal of finite length of $R$.
\end{Theorem}

\begin{pf}
Since $H_{R_+}^i(L) = 0$ for $i > 0$, we get $H_{R_+}^i(R) \cong
H_{R_+}^i(R/L)$ for $i > 0$. Hence $\geom(R) =
\geom(R/L)$. Similarly, since $L+zR/zR$ is an ideal of finite length in
$R/zR$, we also have
$$\geom(R/zR)=\geom((R/zR)/L(R/zR))=\geom((R/L)/z(R/L)).$$ Further it is clear that $p_R(X) = p_{R/L}(X)$ thus,
by replacing $R/L$ with $R,$ we may assume
that $z$ is a regular element in $R$. With this assumption we have
$\geom(R) = \reg(R) = \reg(R/zR)$. We need to show that if $R/zR$ is geometrically $m$-regular,
then $R/zR$ is $(m+p_R(m)-h_R(m))$-regular.

Since $R/zR$ is geometrically $m$-regular, we have  $H^i(R/zR)_{n-i+1}=0$ for every $n 
\ge m$, $i>0$. If we can prove $H_{R_+}^0(R/zR)_{n+1} = 0$  for some 
$n \ge m$, then $R/zR$ is weakly $n$-regular, hence $n$-regular by 
Corollary \ref{basic1}.
So it is sufficient to show that if $R/zR$ is geometrically $m$-regular, there exists an integer $j$ with $m 
+1 \le j \le m+p_R(m)-h_R(m)+1$ such that 
$H_{R_+}^0(R/zR)_j = 0$. \par

From the short exact sequence
$0 \To  R(-1) \overset z \To  R \To  R/zR \To  0$
we get a long exact sequence of local cohomology
modules:
\begin{eqnarray*} 0 \To H_{R_+}^0(R/zR)_n \To H_{R_+}^1(R)_{n-1} \To
H_{R_+}^1(R)_n  \To H_{R_+}^1(R/zR)_n \To \cdots \\ \To
H_{R_+}^i(R/zR)_n \To  H_{R_+}^{i+1}(R)_{n-1} \To H_{R_+}^{i+1}(R)_n \To
\cdots \end{eqnarray*}
for every integer $n$. Since $R/zR$ is geometrically $m$-regular, we have 
$H_{R_+}^i(R/zR)_n = 0$ for ever $i> 0$ and $n \ge m-i+1$. Putting this into the above 
exact sequence yields, with the assumption $n \ge m$, 
\begin{eqnarray}\label{c} \dim_kH_{R_+}^1(R)_{n-1} -\dim_kH_{R_+}^1(R)_n = \dim_kH_{R_+}^0(R/zR)_n
\end{eqnarray}
and the injectivity of the map
$H_{R_+}^{i+1}(R)_{n-1} \To H_{R_+}^{i+1}(R)_n$ for $i \ge 1$.
 Because 
$H_{R_+}^{i+1}(R)_t = 0$ for $t \gg 0$, it follows that  $H_{R_+}^{i+1}(R)_t = 0$ for  $i \ge 1$ and $t\ge m-1$.\par

Since $H_{R_+}^0(R) = 0$, this implies $$h_R(m)-p_R(m)=\sum_{i\ge 0}(-1)^i\dim_kH_{R_+}^i(R)_m=
-\dim_kH_{R_+}^1(R)_m.$$

Put $s = 1+m+p_R(m)-h_R(m)=1+m+\dim_kH_{R_+}^1(R)_m$ and assume by contradiction that 
$H_{R_+}^0(R/zR)_j \neq 0$ for every integer $j$ with $m+1 \le j \le 
s$. By (\ref{c}) we get 

$$\dim_kH_{R_+}^1(R)_m =
\sum_{j=m+1}^s\dim_kH_{R_+}^0(R/zR)_j+\dim_kH_{R_+}^1(R)_s\ge s-m>  \dim_kH_{R_+}^1(R)_m.$$

 The proof of Theorem \ref{Mumford} is now complete.
\end{pf}

\section{Extended degree versus  geometric regularity}

Let $(A,\mm)$ be a local ring with infinite residue field.
Let $\M(A)$ denote the class of finitely generated $A$-modules.
Following [DGV] and [V2], an {\it extended degree} (or cohomological
degree) on $\M(A)$
is a numerical function $D(\cdot)$ on $\M(A)$ such that the
following properties hold for every module
$M \in \M(A)$:\par

(i) $D(M) = D(M/L) + \ell(L)$, where $L$ is the maximal submodule of
$M$ having finite length, \par

(ii) $D(M) \ge D(M/xM)$ for a generic element $x$ of $\mm$,\par

(iii) $D(M) = e(M)$ if $M$ is a Cohen-Macaulay $A$-module, where
$e(M)$ denotes the multiplicity of $M$ with respect to $\mm$. \sk

\noindent{\bf Example.} The prototype of a extended degree is the
{\it homological degree} $\hdeg(M)$ introduced and studied by
Vasconcelos in [V1] (see also [V2]).\par
If $A$ is a homomorphic image of a Gorenstein ring $S$ with $\dim S = 
n$ and $M \in \M(A)$ with $\dim M = r$, we define
$$\hdeg(M) := e(M) + \displaystyle \sum_{i=0}^{r-1} {r-1 \choose
i}\hdeg(\Ext_S^{n-i}(M,S)).$$
This is a recursive definition on the dimension since $\dim 
\Ext_S^{n-i}(M,S) < r$ for $i = 0,\ldots,r-1$.\par

If $A$ is not a homomorphic image of a Gorenstein ring, we only need to  put
$$\hdeg(M) := \hdeg(M \otimes_A\hat A),$$
where $\hat A$ denotes the $\mm$-adic completion of $A$.\par
In particular,
if $M$ is a generalized Cohen-Macaulay module, that is,
$\ell(H_\mm^i(M)) < \infty$
for $i < r = \dim M$, where $H_\mm^i(M)$ denotes the $i$-th local
cohomology module of $M$ with support $\mm$, then
$$\hdeg(M) = e(M)+ \sum_{i=0}^{r-1}{r-1 \choose i}\ell(H_\mm^i(M)).$$
This class of modules is rather large. In fact, if $A$ is a quotient
ring of a Cohen-Macaulay ring, then $M$ is a generalized
Cohen-Macaulay module if and only if  $M$ is locally Cohen-Macaulay
on the punctured spectrum of $A$ and $\operatorname{Supp}(M)$ is
equidimensional. \sk

Any extended degree $ D(M) $ will satisfy $ D(M) \ge e(M) $ with
equality holding if and only if $M$ is a
Cohen-Macaulay module. Following [DGV]  we call  the difference
$$I(M) : = D(M)-e(M).$$
a {\it Cohen-Macaulay deviation} of $M$. It is obvious that $I(M)$
satisfies the following conditions: \par
(i') $I(M) = I(M/L) + \ell(L)$,  \par
(ii'') $I(M) \ge I(M/xM)$ for a generic element $x$ of $\mm$.\sk

The Hilbert-Samuel function can be bounded in terms of any extended
degree as follows.  \sk

\begin{Theorem} \label{uniform}
Assume that $d = \dim A \ge 1$. Let  $I(A)$ be an arbitrary
Cohen-Macaulay deviation of $A$. For all $n \ge 0$ we have
$$\ell(A/\mm^{n+1}) \le e(A){n+d-1 \choose d} +  I(A){n+d-2 \choose
d-1} + {n+d-1 \choose d-1}.$$
\end{Theorem}

\begin{pf} Let $x_1,\ldots,x_d$ be a system of generic elements  of
$\mm$ such that $e(A)$ is the multiplicity of $A$ with respect to the
ideal $\qq = (x_1,\ldots,x_d)$. We have
$$\ell(A/\mm^{n+1}) \le \ell(A/\mm\qq^n) = \ell(A/\qq^n) +
\ell(\qq^n/\mm\qq^n).$$
Since $x_1,\ldots,x_d$ are analytically independent,
$$\ell(\qq^n/\mm\qq^n)  =  {n+d-1 \choose d-1}.$$
It remains to show that
$$\ell(A/\qq^n) \le e(A) {n+d-1\choose d} + I(A){n+d-2 \choose d-1}.$$\par

If $d = 1$, $A/L$ is a Cohen-Macaulay ring, where $L$ is the largest
ideal of finite length.  Hence $\ell(A/\qq^n+L) = e(A)n$ and $I(A) =
\ell(L)$.
  From this it follows that
$$\ell(A/\qq^n) \le \ell(A/\qq^n+L) + \ell(L) = e(A)n + I(A).$$\par

If $d > 1$, we put $\bar A = A/(x_1)$ and $\bar \qq = \qq/(x_1)$.
Then $\dim \bar A = d-1 \ge 1$, $e(\bar A) = e(A)$ and $I(\bar A) \le
I(A)$. By induction we may assume that
$$\ell(\bar A/\bar \qq^i) \le e(A){i+d-2 \choose d-1} + I(A){i+d-3
\choose d-2}.$$
for all $i \ge 1$. From the exact sequence
$$0 \To \qq^{i-1}:x_1/\qq^{i-1}  \To A/\qq^{i-1} \overset {x_1}  \To
A/\qq^i \To \bar A/\bar \qq^i \To 0$$
we can deduce that
$$\ell(\qq^{i-1}/\qq^i) = \ell(A/\qq^i)-\ell(A/\qq^{i-1}) \le
\ell(\bar A/\bar \qq^i).$$
Using these formulas for $i = 1,\ldots,n-1$ we get
\begin{align*}
\ell(A/\qq^n) & = \sum_{i=1}^n \ell(\qq^{i-1}/\qq^i) \le
\sum_{i=1}^{n-1} \ell(\bar A/\bar \qq^i)\\
& \le \sum_{i=1}^n \big[e(A){i+d-2 \choose d-1} + I(A){i+d-3 \choose
d-2}\big]\\
& = e(A){n+d-1 \choose d} + I(A){n+d-2 \choose d-1}.
\end{align*}
\end{pf}

\noindent{\bf Remark.} From earlier results on the Hilbert functions of local
rings
[Tru1, Corollary 2.2], [DGV, Theorem 4.6], [RVV, Theorem 2.2] one can
only derive the bound:
$$\ell(A/\mm^{n+1}) \le D(A){n+d-1 \choose d} + {n+d-1 \choose d-1}.$$
This bound is much weaker than the bound of Theorem \ref{uniform}.\sk

We will use the above uniform bound for the Hilbert-Samuel function
to give a local version of Mumford's Theorem \ref{Mumford}. For this
we shall need the following observations.  \sk

Let $G$ denote the associated graded ring $\bigoplus_{n\ge
0}\mm^n/\mm^{n+1}$ of $A$. Let $x$ be a generic element of $\mm$.
Without restriction we may assume that the initial form $x^*$ of $x$
in $G$ is a filter-regular element.
Let $\bar G$ denote the associated graded ring of  $A/(x)$.   \sk

\begin{Lemma} \label{superficial} {\rm (cf. [ST2, Lemma 1])}
$\geom(G/(x^*)) = \geom(\bar G).$  \end{Lemma}

\begin{pf} We have
\begin{eqnarray*} & G/(x^*)  =  \bigoplus_{n\ge 0}
\mm^n/(\mm^{n+1}+ x\mm^{n-1}), &\\ &
\bar G = \bigoplus_{n\ge 0} (\mm^n + (x))/(\mm^{n+1} + (x))=
\bigoplus_{n\ge 0} \mm^n/(\mm^{n+1} + (x) \cap \mm^n).
\end{eqnarray*} Hence there is a natural graded epimorphism  from
$G/(x^*)$ to $\bar G$ whose kernel is
$$E = \bigoplus_{n \ge 0}(\mm^{n+1} + (x)\cap \mm^n)/(\mm^{n+1} +
x\mm^{n-1}).$$
Since $x$ is a superficial element, $\mm^{n+1}:x = \mm^n + (0:x)$ for
all large $n$  [N, Theorem (22.6)]. Thus,  $(x) \cap \mm^{n+1} =
x\mm^n$ so that $E_n = 0$ for all large $n$. This implies
$H_{G_+}^i(G/(x^*)) =
H_{G_+}^i(\bar G)$ for $i > 0$. Hence $\geom(G/(x^*)) =
\geom(\bar G).$ \end{pf}

For $n \ge \reg(\bar G)$, we can bound $p_G(n)$ by means of the
Hilbert-Samuel function of $A/(x)$. Note that $\reg(\bar G) \ge
\geom(\bar G)$.\sk

\begin{Lemma} \label{reduction}
$p_G(n) = \ell(A/\mm^{n+1}+(x))-\ell(0:x)$ for $n \ge \reg(\bar G)$.
\end{Lemma}

\begin{pf} By [N, Theorem (22.6)] we have
$$h_G(n) = \ell(\mm^n/\mm^{n+1}) = \ell(A/\mm^{n+1}+(x)) - \ell(0:x)$$
for all large $n$. Write $\ell(A/\mm^{n+1}+(x)) = \sum_{i = 0}^nh_{\bar G}(i).$
By Lemma \ref{Serre}, $h_{\bar G}(i) = p_{\bar G}(i)$ for $i >
\reg(\bar G)$. Hence $\ell(A/\mm^{n+1}+(x))-\ell(0:x)$ is a
polynomial function for $n \ge \reg(\bar G)$. This polynomial must be
$p_G(n)$. \end{pf}\sk

Now we can deduce from Mumford's Theorem \ref{Mumford} the following
local version.

\begin{Theorem} \label{local}
Assume that $d = \dim A \ge 2$. Let  $I(A)$ be an arbitrary
Cohen-Macaulay deviation of $A$. For $n \ge \reg(\bar G)$ we have
$$\geom(G) \le e(A){n+d-2 \choose d-1} + I(A){n+d-3 \choose d-2}.$$
\end{Theorem}

\begin{pf} By Theorem \ref{Mumford} and Lemma \ref{superficial} we have
$$\geom(G) \le n + p_G(n) - h_{G/L}(n),$$
where $L$ is the largest ideal of finite length of $G$. Since $\dim 
G/L = d$,  $$h_{G/L}(n) \ge {n+d-1 \choose d-1}.$$
By Lemma \ref{reduction} and Theorem \ref{uniform},
\begin{align*}
  p_G(n) &  \le  \ell(A/\mm^{n+1}+(x))\\
& \le e(A/(x)){n+d-2 \choose d-1} +
I(A/(x)){n+d-3 \choose d-2} + {n+d-2 \choose d-2}.
\end{align*}
Since $e(A/(x)) = e(A)$ and $I(A/(x)) \le I(A)$, we finally obtain
\begin{align*}
\geom(G) &  \le  n + e(A){n+d-2 \choose d-1} +  I(A){n+d-3 \choose 
d-2}+ {n+d-2 \choose d-2} - {n+d-1 \choose
d-1} \\ & \le e(A){n+d-2 \choose d-1} + I(A){n+d-3 \choose d-2}.
\end{align*}
\end{pf}

\section{Bounds for the regularity}

Let $(A,\mm)$ be a local ring with infinite residue field and $d =
\dim A$. Let $G$ be the associated graded ring of $A$. We will apply
Theorem \ref{local} to compute $\reg(G)$. For that we shall need the
following two propositions.  The first proposition allows us to pass
to the case $\depth A > 0$, while the second proposition shows that
regularity coincides with geometric regularity in this case. \sk

\begin{Proposition} \label{depthzero}
Let $G'$ denote the associated graded ring of $A/L$, where $L$ is the
largest ideal of $A$ of finite length.  Then
$$\reg(G) \le \reg(G') + \ell(L).$$
\end{Proposition}

\begin{pf} We have
$$G' = \bigoplus_{n\ge 0} (\mm^n+L)/(\mm^{n+1}+L) \cong \bigoplus_{n\ge
0}\mm^n/(\mm^{n+1} + \mm^n\cap L).$$
Therefore, there is a natural graded epimorphism from $G $ to
$G'$ with the kernel
$$K = \bigoplus_{n \ge 0}(\mm^{n+1} + \mm^n\cap L)/\mm^{n+1}
\cong \bigoplus_{n\ge 0}\mm^n\cap L/\mm^{n+1}\cap L.$$
Note that
$$\ell(K) = \sum_{n \ge 0} \ell(\mm^n\cap L/\mm^{n+1}\cap L) = \ell(L).$$
Then $K$ has finite length. This implies $H_{G_+}^i(G) \cong
H_{G_+}^i(G')$ for $i >  0.$\par

Put $a = \reg(G')$ and $\ell = \ell(K)$. Then there exists an index $m$
with $a \le m \le a + \ell$ such that
$K_{m+1} = 0.$ Since $G'$ is $m$-regular, $H_{G_+}^i(G)_{m-i+1}=
H_{G_+}^i(G')_{m-i+1} = 0$ for
$i \ge 0$. Hence $G$ is weakly
$m$-regular. By Corollary \ref{basic1}, $G$ is
$m$-regular so that $\reg (G)\le m \le a + \ell.$
\end{pf} \sk

\begin{Proposition} \label{Hoa}
Assume that $\depth A > 0$. Then
$$\reg(G) = \geom(G).$$
\end{Proposition}

\begin{pf} Note that the result is trivial if $\depth G > 0.$
For any integer $i$ we set $a_i(G) :=
\max\{n|\ H_{G_+}^i(G)_n
\neq 0\},$ where $a_i(M) = -\infty$ if
$H_{G_+}^i(G) = 0$. Then
\begin{eqnarray*} \reg(G) & = & \max\{a_i(G)+i|\ i \ge 0\},\\ \geom(G) & =
& \max\{a_i(G)+i|\ i > 0\}. \end{eqnarray*}
By [H, Theorem 5.2] (see also [Ma, Theorem 2.1]),
the assumption $\depth A > 0$ implies
$a_0(G) \le a_1(G)$. Hence
$\reg(G) = \geom(G).$  \end{pf} \sk

Now we are able to give an upper bound for the Castelnuovo-Mumford regularity
of the associated graded ring in terms of any extended degree.
\sk

\begin{Theorem} \label{regularity} Let $I(A)$ be an arbitrary
Cohen-Macaulay deviation of $A$. Then\par
{\rm (i) } $\reg(G) \le e(A) + I(A)-1$ if $d = 1$,\par
{\rm (ii)} $\reg(G) \le
   e(A)^{(d-1)!-1}[e(A)^2+e(A)I(A)+2I(A)-e(A)]^{(d-1)!}- I(A)$ if $d \ge 2$.
\end{Theorem}

\begin{pf} Let $G'$  and $L$ be as in Proposition \ref{depthzero}.
Then $\reg(G) \le \reg(G') + \ell(L).$ Since $e(A) = e(A/L)$ and
$I(A) = I(A/L) + \ell(L)$, we only need prove the conclusion for the
ring $A/L$. Replacing $A$ by $A/L$ we may assume that $\depth A > 0$.
\par

If $d = 1$, then $A$ is a Cohen-Macaulay ring. Hence $I(A) = 0$.  By
[Tru2,  Theorem 5.1(i) and Theorem 1.2(iii)]  we already know that
$\reg(G) \le e(A)-1$ (see also the proof of [ST1, Lemma 5]).  \par

For $d \ge 2$ we choose a generic element $x$ of $\mm$ such that the
initial form $x^*$ of $x$ is a filter-regular linear form in $G$.
Note that $e(A/(x)) = e(A)$ and $I(A/(x)) \le I(A)$. Let $\bar G$
denote the associated graded ring of $A/xA$.
Put $m = \reg(\bar G)$. By Proposition \ref{Hoa} and Theorem
\ref{local},  we have
$$\reg(G) = \geom(G) \le e(A){m+d-2 \choose d-1} + I(A){m+d-3 \choose d-2}.$$
Since
$$ {m+d-2\choose d-1}  \le  m^{d-1},\ \ \  {m+d-3 \choose d-2}  \le  m^{d-2},$$
this implies
$$\reg(G) \le e(A)m^{d-1} + I(A)m^{d-2}.$$ \par

If $d = 2$, we have
$$m \le e(A/xA) + I(A/xA)-1 \le  e(A)+I(A)-1.$$
  From this it follows
\begin{align*}
\reg(G)  \le  e(A)m + I(A) & = e(A)[e(A) + I(A)-1] + I(A)\\
& =  [e(A)^2+e(A)I(A)+ 2I(A) - e(A)]-I(A).
\end{align*}  \par

If $d \ge 3$, we have
$$ e(A)m^{d-1} +  I(A)m^{d-2} \le  e(A)[m+I(A)]^{d-1}- I(A). $$
Using induction we may assume that
$$m \le e(A)^{(d-2)!-1}[e(A)^2+e(A)I(A)+2I(A)-e(A)]^{(d-2)!}- I(A) .$$
Therefore, the last bound for $\reg(G)$ implies
\begin{align*}
\reg(G) & \le
e(A)\big\{e(A)^{(d-2)!-1}[e(A)^2+e(A)I(A)+2I(A)-e(A)]^{(d-2)!}\big\}^{d-1}-
I(A)
\\ & \le  e(A)^{(d-1)!-1}[e(A)^2+e(A)I(A)+2I(A)-e(A)]^{(d-1)!}- I(A).
\end{align*}
The proof of Theorem \ref{regularity} is now complete. \end{pf}
\sk

\begin{Corollary} Let $A$ be a Cohen-Macaulay local ring with $d =
\dim A \ge 1$. Then \par
{\rm (i) } $\reg(G) \le e(A) -1$ if $d = 1$,\par
{\rm (ii)}
   $\reg(G) \le e(A)^{2((d-1)!)-1}[e(A)-1]^{(d-1)!}$ if $ d \ge 2$.
\end{Corollary}

\begin{pf} The Cohen-Macaulayness of $A$ implies $I(A) = 0$. Hence
the conclusion follows from Theorem \ref{regularity}. \end{pf} \sk

As one can see from the proof, the bound of Theorem \ref{regularity}
can be further improved if $d \ge 3$, but the formula is not so
compact.

\section{Finiteness of Hilbert functions}

Let $(A,\mm)$ be a local ring with infinite residue field and $d =
\dim A$.  We can bound the coefficients $e_i(A)$ in
terms of any extended degree as follows.

\begin{Theorem} \label{coeff}
Let $I(A)$ be an arbitrary Cohen-Macaulay deviation of $A$. Then \par
{\rm (i)} $|e_1(A)| \le  \displaystyle \frac{e(A)[e(A)-1]}{2}+I(A)$,\par
{\rm (ii)} $|e_i(A)| \le e(A)^{i!-i}[e(A)^2+e(A)I(A)+2I(A)]^{i!}-1$ if $i
\ge 2$.
\end{Theorem}

\begin{pf}
Let $L$ be the largest ideal of finite length of $A$. Note that
\begin{eqnarray*}
\ell(A/\mm^{n+1}) & = & \ell(A/\mm^{n+1}+L) +
\ell(\mm^{n+1}+L/\mm^{n+1})\\ & = &
\ell(A/\mm^{n+1}+L) + \ell(L/\mm^{n+1}\cap L)\ =\ \ell(A/\mm^{n+1}+L) +
\ell(L)\end{eqnarray*}
for $n$ large. Then we have
\begin{eqnarray*}
e_i(A) & = & e_i(A/L),\ i = 0,\ldots,d-1,\\
|e_d(A)| & \le & |e_d(A/L)|  + \ell(L).
\end{eqnarray*}\par

If $d = 1$, then $A/L$ is a Cohen-Macaulay ring. It is easy to see (see e.g
[Ki]) that
$$|e_1(A/L)| \le
\frac{e(A/L)[e(A/L)-1]}{2}.$$
Since $e(A) = e(A/L)$ and $I(A) = \ell(L)$, we get
$$|e_1(A)| \le |e_1(A/L)| + \ell(L) \le \frac{e(A)[e(A)-1]}{2} + I(A).$$
\par

If $d \ge 2$, we first consider the case $\depth A > 0$. Choose a generic
element $x$ of $\mm$ such that the initial form $z$ of $x$ is a
filter-regular linear form in $G$. Since $x$ is a superficial element of
$\mm$, we have $\dim A/xA
= d-1$, $e_i(A) = e_i(A/x)$ for $i = 0,...,d-1$ (see e.g. [N, Theorem 
(22.6)] or [ST2, Lemma 1]) and
$I(A/(x)) \le I(A)$. Using the induction hypothesis on $A/(x)$ we may
assume that
\begin{eqnarray*} |e_1(A)| & \le & \frac{e(A)[e(A)-1]}{2} + I(A),\\
|e_i(A)| & \le &
   e(A)^{i!-i}[e(A)^2+e(A)I(A)+2I(A)]^{i!}-1, \  i = 2,...,d-1.
\end{eqnarray*} \par

It remains to prove the bound for $e_d(A)$. We have
$$(-1)^de_d(A) = P_A(n)-\sum_{i=0}^{d-1}(-1)^ie_i(A){m+d-i \choose d-i}$$
for all $n > 0$.
Let $G$ be the associated graded ring of $A$. Since  $\ell(A/\mm^{n+1}) =
\sum_{i=0}^nh_G(i)$, Lemma \ref{Serre} yields $P_A(n) = \ell(A/\mm^{n+1})$
for $n \ge \reg(G)$. Put
$$m = e(A)^{(d-1)!-1}[e(A)^2+e(A)I(A)+2I(A)]^{(d-1)!}-1.$$
Then $m \ge \reg(G)$ by Theorem \ref{regularity}.
Therefore, using Theorem \ref{uniform} we have
$$|e_d(A)|  \le I(A){m+d-2 \choose d-1} + {m+d-1 \choose d-1}
+  {m+d-1 \choose d }e - {m+d \choose d }  e +$$ $$+
\sum_{i=1}^{d-1}|e_i(A)|{m+d-i  \choose d-i}
    \le I(A){m+d-2 \choose d-1} + {m+d-1 \choose d-1}
+   \sum_{i=1}^{d-1}|e_i(A)|{m+d-i  \choose d-i}.$$
Note that
$${m+d-2 \choose d-1} \le m^{d-1},\ \  \  {m + d-i  \choose d-i}  \le
(d-i+1) m^{d-i}. $$
A numerical computation shows that
\begin{align*}
|e_1(A)| & \le \frac{e(A)[e(A)-1]}{2} + I(A) \le \frac{m -
I(A)-d+1}{d}.\end{align*}
Moreover
\begin{align*}
|e_i(A)| & \le m,\ i = 2,\ldots,d-1.
\end{align*}
Then we get
\begin{align*}
|e_d(A)| & \le  I(A)m^{d-1} + dm^{d-1} +   \frac{m-I(A)-d+1}{d} \cdot
d m^{d-1} +
\sum_{i=2}^{d-1}m^{d-i+1}(d-i+1)\\
& = m^d + dm^{d-1} + (d-2)m^{d-2} + \cdots + 2m^2 \le   (m+1)^d - 1\\
&  = e(A)^{d!-d}[e(A)^2+e(A)I(A)+2I(A)]^{d!}-1.
\end{align*}\par

Finally, we consider the case $\depth A = 0$. As shown above, we have
\begin{align*}
|e_1(A/L)|  & \le  \displaystyle \frac{e(A/L)[e(A/L)-1]}{2}+I(A/L),\\
|e_i(A/L)|  & \le e(A/L)^{i!-i}[e(A/L)^2+e(A/L)I(A/L)+2I(A/L)]^{i!}-1,\ i =
2,...,d.
\end{align*}
Since $e(A) = e(A/L)$ and $I(A) = I(A/L) + \ell(L)$, this immediately
implies
$$
|e_1(A)| =  |e_1(A/L)|  \le  \displaystyle \frac{e(A)[e(A)-1]}{2}+I(A),$$
$$|e_i(A)|  =  |e_i(A/L)|  \le
e(A )^{i!-i}[e(A)^2+e(A)I(A)+2I(A)]^{i!}-1,\ i = 2,...,d-1, $$
$$|e_d(A)|  \le |e_d(A/L)|  + \ell(L) \le
e(A)^{d!-d}[e(A)^2+e(A)I(A/L)+2I(A/L)]^{d!}-1 + \ell(L) \le $$
$$ e(A)^{d!-d}[e(A)^2+e(A)I(A )+2I(A)]^{d!}-1.$$
   \par
The proof of Theorem \ref{coeff} is now complete.
\end{pf}

\begin{Corollary} \label{Macaulay}
Let $A$ be a Cohen-Macaulay ring. Then
\par {\rm (i)} $|e_1(A)| \le  \displaystyle \frac{e(A)[e(A)-1]}{2}$
{\rm (see e.g. [Ki])}, \par
{\rm (ii)} $|e_i(A)| \le e(A)^{3(i!)-i}-1$ if $i \ge 2$.
\end{Corollary}

\begin{pf} This follows from the fact that $I(A) = 0$ for a
Cohen-Macaulay ring $A$. \end{pf}

\noindent{\bf Remark.} If $A$ is  a Cohen-Macaulay ring, Srinivas
and Trivedi
[ST2, Theorem 1] already gave the bound:
$$|e_i(A)| \le (9e(A)^5)^{i!},\; i = 1,...,d.$$
If $A$ is a generalized Cohen-Macaulay ring, this bound has been
extended by Trivedi [Tri2, Theorem 8] to
$$|e_i(A)| \le [(3+c)^2e(A)^5]^{i!}, i = 1,\ldots,d,$$
for some invariant $c \ge 2I(A)$, where $I(A) = \sum_{i=0}^{d-1}{d-1
\choose i}\ell(H_\mm^i(A))$ is the Cohen-Macaulay deviation of the
homological degree $\hdeg(A)$. These bounds are worse than the bounds
of Theorem \ref{coeff} and Corollary \ref{Macaulay}. \sk

As an application of our bounds for the Castelnuovo-Mumford
regularity and the coefficients of the Hilbert-Samuel polynomial we
obtain the finiteness of Hilbert-Samuel functions of local rings with
given dimension and extended degree. \sk

\begin{Corollary} \label{finite} Given two positive integers $d$ and $q$
there exist only a finite number of
Hilbert-Samuel functions for a local ring $A$ with $\dim A=d$ and
$D(A)\le q.$ \end{Corollary}

\begin{pf} By Lemma \ref {Serre} it follows that  $P_A(n) =
\ell(A/\mm^{n+1})$ for $n \ge\reg(G)$. By Theorem
\ref{uniform}, there are only a finite number of possibilities for
$\ell(A/\mm^{n+1})$ for  a fixed $n$.
Hence the finiteness of the number of the possibilities for the
function $\ell(A/\mm^{n+1})$ follows from the
finiteness of possibilities for
$\reg(G)$ and for the polynomial $P_A(n)$. Hence the conclusion
follows from Theorem \ref{regularity} and
Theorem \ref{coeff}. \end{pf}

Obviously, Corollary \ref{finite} implies the finiteness of the
number of Hilbert-Samuel functions of Cohen-Macaulay rings with given
dimension and multiplicity [ST2] and of generalized Cohen-Macaulay
rings with given dimension, multiplicity  and lengths of local
cohomology modules [Tri2]. Moreover, our results
can be modified to cover
the finiteness of the number of Hilbert-Samuel functions of graded
algebras with given dimension and extended
degree [RVV].
\sk

To illustrate the above results we consider the following simple example.\sk

\noindent {\bf Example.}  Let $r \ge 1$ be any integer and consider the
one-dimensional
local ring $A = k[[x,y]]/(x^2,xy^r)$. Then $A$ is a
non-Cohen-Macaulay local ring with multiplicity $e(A) =
1$ and Hilbert-Samuel function
$$\ell(A/\mm^{n+1}) = \left\{\begin{array}{lll} 2n+1 & \text{for} & n \le r,\\
n + r + 1 & \text{for} & n > r. \end{array} \right. $$
This shows that there may be
infinite Hilbert-Samuel functions for local rings with given dimension and
multiplicity (see [ST1, Section 4] for a similar example with
two-dimensional local domains).  Since $L = (x)/(x^2,xy^r)$, $A$ is a
generalized Cohen-Macaulay ring with $I(A) = \ell(L) = r$. Hence
$D(A) = r+1$. \sk

\section*{References}\sk

\noindent [BM] D.~Bayer and D.~Mumford, What can be computed in algebraic
geometry~? in: D. Eisenbud and L. Robbiano
(eds.), Computational Algebraic Geometry and Commutative Algebra,
Proceedings, Cortona (1991), Cambridge University
Press, 1993, 1--48. \par

\noindent [BH] W. Bruns and J. Herzog, Cohen-Macaulay rings, Cambridge,
1998.\par

\noindent [DGV] L.~R.~Doering, T.~Gunston and W.~Vasconcelos, Cohomological
degrees and Hilbert functions of graded modules,
    Amer. J. Math. 120 (1998), 493--504. \par

\noindent [E] D.~Eisenbud, Commutative Algebra with a viewpoint toward
Algebraic Geometry, Springer, 1994. \par

\noindent [EG] D.~Eisenbud and S.~Goto, Linear free resolutions and minimal
multiplicities, J.~Algebra~88 (1984),
89-133.\par

\noindent [Ki] D.~Kirby, The reduction number of a one-dimensional local ring,
J. London Math. Soc. 10 (1975), 471-481.\par

\noindent [Kl] S. Kleiman, Th\'eorie des intersections et th\'eor\`eme de
Riemann-Roch, in: SGA 6, Lect. Notes in
Math. 225, Springer, 1971.\par

\noindent [H] L.T.~Hoa, Reduction numbers of equimultiple ideals, J. Pure
Appl. Algebra 109 (1996), 111-126. \par

\noindent [Ma] T.~Marley, The reduction number of an ideal and the local
cohomology of the associated graded ring,
Proc. Amer. Math. Soc. 117 (1993), 335-341. \par

\noindent [Mu] D.~Mumford, Lectures on curves on an algebraic surfaces,
Princeton Univ. Press, Princeton, 1966.\par

\noindent [N]  M. Nagata, Local rings, Interscience, New York, 1962.\par

\noindent [RVV] M.~E.~Rossi, G.~Valla, and W.~Vasconcelos, Maximal Hilbert
functions, Results in Math. 39 (2001) 99-114.\par

\noindent [ST1] V. Srinivas and V. Trivedi, A finiteness theorems for the
Hilbert functions of complete intersection
local rings, Math. Z. 225 (1997), 543-558. \par

\noindent [ST2] V. Srinivas and V. Trivedi, On the Hilbert function of a
Cohen-Macaulay ring, J. Algebraic Geom. 6
(1997), 733-751. \par

\noindent [Tri1] V. Trivedi, Hilbert functions, Castelnuovo-Mumford regularity
and uniform Artin-Rees numbers,
Manuscripta Math. 94 (1997), no. 4, 485--499. \par

\noindent [Tri2] V. Trivedi, Finiteness of Hilbert functions for
generalized Cohen-Macaulay modules, Comm. Algebra 29 (2) (2001),
805--813. \par

\noindent [Tru1] N.~V.~Trung, Absolutely superficial sequence, Math. Proc.
Camb. Phil. Soc. 93 (1983), 35-47. \par

\noindent [Tru2] N.~V.~Trung, Reduction exponent and degree bound
for the defining equations of graded rings, Proc. Amer. Math.
Soc.  101 (1987), 229-236. \par

\noindent [V1] W. Vasconcelos, The homological degree of a module, Trans.
Amer. Math. Soc. 350 (1998), no. 3,
1167--1179.\par

\noindent [V2] W. Vasconcelos, Cohomological degrees of graded modules. Six
lectures on commutative algebra
(Bellaterra, 1996), 345--392, Progr. Math. 166, Birkh\"auser, Basel, 1998. \par

\end{document}